\documentclass[a4paper,twoside]{article}
\usepackage[utf8]{inputenc}
\usepackage{float}
\usepackage{units}
\usepackage{textcomp}
\usepackage{amsmath}
\usepackage{amsthm}
\usepackage{amssymb}
\usepackage{graphicx}
\usepackage[numbers]{natbib}

\makeatletter

\pdfpageheight\paperheight
\pdfpagewidth\paperwidth

\numberwithin{equation}{section}
\numberwithin{figure}{section}
\newcommand{\lyxaddress}[1]{
\par {\raggedright #1
\vspace{1.4em}
\noindent\par}
}

\@ifundefined{date}{}{\date{}}

\usepackage{latexsym}
\usepackage{a4wide}
\pdfoutput=1

\RequirePackage{fix-cm}

\usepackage[T1]{fontenc}
\usepackage{array}
\usepackage{url}
\usepackage{multirow}
\usepackage{amsmath}
\usepackage{booktabs}
\usepackage{graphicx}

\renewcommand{\rho}{\varrho}
\renewcommand{\phi}{\varphi}

\theoremstyle{remark}

\makeatother

\begin{document}

\title{Similarity Dimension of Fractal Curves with Multiple Generators}

\author{Stefan Pautze}

\maketitle

\lyxaddress{\begin{center}
Visualien der Breitbandkatze, Pörnbach, Germany, stefan@pautze.de
\par\end{center}}
\begin{abstract}
We propose a definition for the similarity dimension of fractal curves
with multiple generators.
\end{abstract}

\section{Introduction}

Most fractal curves as the Koch curve \citep{HvK1904,HvK1906,Mandelbrot77a}
are defined by an initiator and one generator. While the initiator
is defined by a polygonal chain or just a single line segment, the
generator is defined by a polygonal chain with at least two line segments.
The fractal curve is generated iteratively. In a first step all line
segments of the initiator are replaced by scaled down copies of the
generator with defined orientations. In the following steps, hence
iterations, the line segments are replaced again accordingly. To obtain
a continuous curve that cannot be differentiated anywhere, the process
is repeated an infinite number of times.

In literature a large number of fractal curves can be found, but all
fractal curves known to the author rely on one generator only. An
exception is the FASS (space filling, self avoiding, self similar
and simple) curve of the regular pentagon recently derived in \citep{Pautze2021}.
In detail appropriate decorations were applied on all prototiles of
an appropriate cyclotomic aperiodic substitution tiling (CAST) which
is also a stone inflation.

\begin{figure}[H]
\begin{center}
\resizebox{0.95\textwidth}{!}{%

\includegraphics{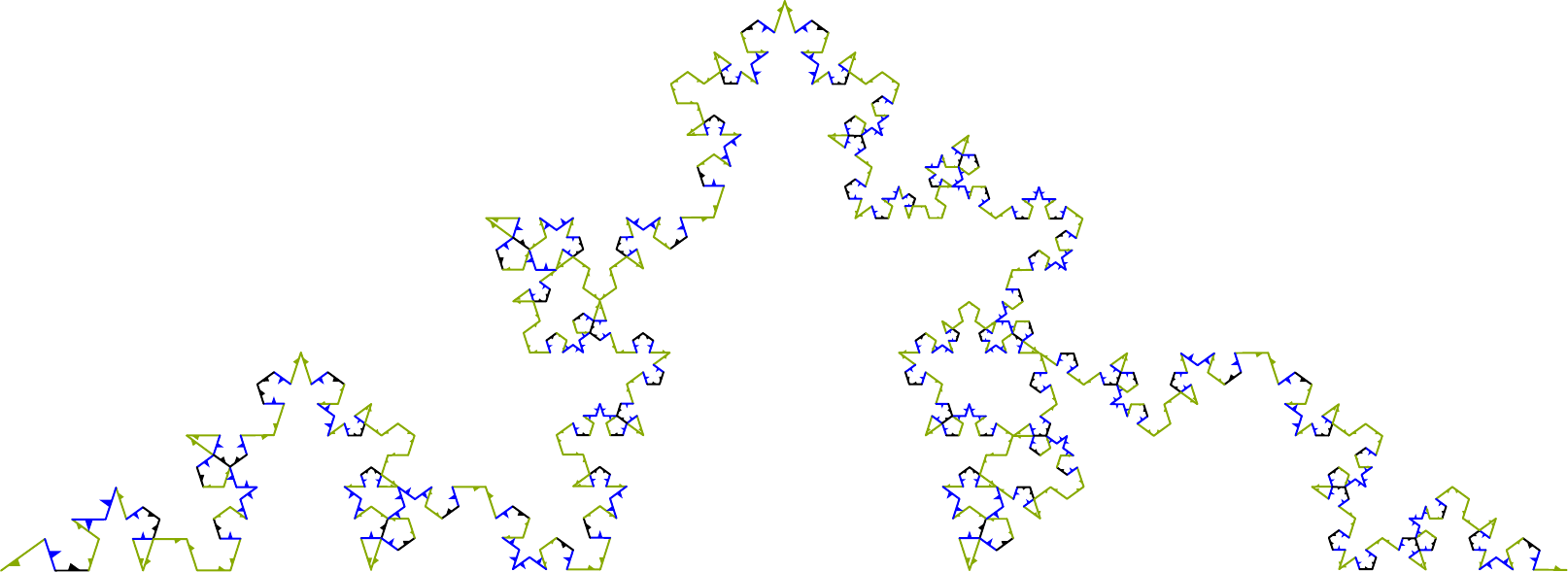}}
\end{center}\caption{\label{curve}The figure shows the $7$th iteration of a fractal curve
with multiple generators as shown in Figure~\ref{generators}. The
different line segment types are marked in green, blue and black.
The initiator of the curve is line segment $L_{1}$.}
\end{figure}

Aperiodic substitution tilings as described in \citep[and references therein]{oro38933,HFonl}
and its cyclotomic variant as discussed in \citep[and references therein]{sym9020019}
are defined by a set of prototiles which can be expanded with a linear
map - the “inflation multiplier” - and dissected into copies of prototiles
in original size - the “substitution rule”. Substitution tilings allow
to cover the entire Euclidean plane without gaps and overlaps with
tiles of finite size or to tile a finite area with tiles of infinitesimal
size. For the latter case just the expansion is omitted. Obviously
in both cases similarity dimension of the tiling is $D=2$. A space
filling curve derived from the tiling has the same property. Since
most substitution tilings rely on $n\geq2$ proto tiles and substitution
rules we can also describe the corresponding FASS curve by a set of
$n$ substitution rules or more precisely $n$ generators. On closer
inspection it turns out that the terms ``iteration'' and ``substitution''
as well as the terms ``generator'' and ``substitution rule'' have
almost identical meanings here.

A FASS curve regardless of the number generators is just a special
case of a fractal curve with similarity dimension $D=2$. This raises
an interesting question: How to derive the similarity dimension of
a fractal curve with multiple generators in general?

In the following sections we will derive a proposal for an appropriate
definition.

\section{Perron-Frobenius Theorem\label{sec:Perron-Frobenius-Theorem}}

For the analysis of aperiodic substitution tilings it is common to
use substitution matrices based on the substitution rules. We will
apply the same principle to fractal curves with multiple generators.
The substitution matrix $M$ describes, how each of the line segments
$L_{i};\,1\leq i\leq n$ is replaced by an individual polygonal chain
$G_{i}$ (the generator), each made of a number $k_{i,1},k_{i,2}\,\ldots\,k_{i,n}$
of line segments $L_{1},L_{2}\,\ldots\,L_{n}$ scaled down by a linear
map with a scaling factor $0<r_{i}<1;\,r_{i}\in\mathbb{R}$. The scaling
factor $r_{i}$ is defined as the ratio of the length of a line segment
$L_{i}$ to the distance between start and end point of the polygonal
chain of generator $G_{i}$.

As long as $M\geq0$ is a primitive matrix so that $\exists\left(k,M\geq0\right)\Rightarrow M^{k}>0$
we can apply the Perron–Frobenius theorem as introduced in \citep{Perron1907}
and \citep{Frobenius1912}. That means the Perron–Frobenius eigenvalue
of the substitution matrix $M$ is a real positive number so that
$\lambda_{PF}\in\mathbb{R}_{>0}$ and the modulo of every other eigenvalue
is smaller than $\lambda_{PF}$. Furthermore the corresponding left
and right Perron-Frobenius eigenvectors $x_{PFL}$ and $x_{PFR}$
also have only real positive elements.

As discussed in \citep[and references therein]{oro38933} the elements
of a normalized right Perron-Frobenius eigenvector $x_{PFR}^{*}$
which is defined as

\begin{equation}
x_{PFR}^{*}=\frac{x_{PFR}}{\left\Vert x_{PFR}\right\Vert _{1}}=\left(\begin{array}{c}
f_{1}\\
f_{2}\\
\vdots\\
f_{n}
\end{array}\right)\ so\ that\ \sum_{i=1}^{n}f_{i}=1\label{eq:NRPFE}
\end{equation}
describe the relative frequency of the proto tiles - or applied to
our problem - the relative frequency of the different types of line
segments $L_{1},L_{2}\,\ldots\,L_{n}$ in a fractal curve. Furthermore
the Perron–Frobenius eigenvalue $\lambda_{PF}$ represents the multiplication
factor applied to the total number of tiles during the substitution.
As a consequence it also represents the multiplication factor applied
to the total number of line segments during an iteration.

It is known that repeatedly multiplying a vector $x\geq0$ by a primitive
matrix $M$ cause the result to converge to the Perron-Frobenius eigenvector
of the matrix. That means regardless of the initial setting the relative
frequencies of the types of lines segments of the fractal curves will
converge with an increasing number of iterations to the relative frequencies
described by $x_{PFR}^{*}$.

\section{Similarity Dimension}

During the research of fractals and fractal structures different types
of fractal dimensions were defined by well known mathematicians like
F. Hausdorff, H. Minkowsky, B. B. Mandelbrot and many more. For strictly
self similar structures as discussed in this article it is sufficient
to focus on the similarity dimension.

H.-O. Peitgen et. al. defined strict self similarity in \citep{Peitgen1991}
as follows: ``If the figure can be decomposed into parts which are
exact replicas of the whole, then the figure is called strictly self-similar.
Every part of a strictly self-similar structure contains an exact
replica of the whole.'' For the purpose of this paper we assume that
strict self similarity also applies to a set of figures which can
be decomposed into a set of parts which are exact replicas of elements
in the set of wholes.

The similarity dimension is defined as follows:

\begin{equation}
1=nr^{D}
\end{equation}

\begin{figure}
\begin{center}
\resizebox{0.85\textwidth}{!}{%

\includegraphics{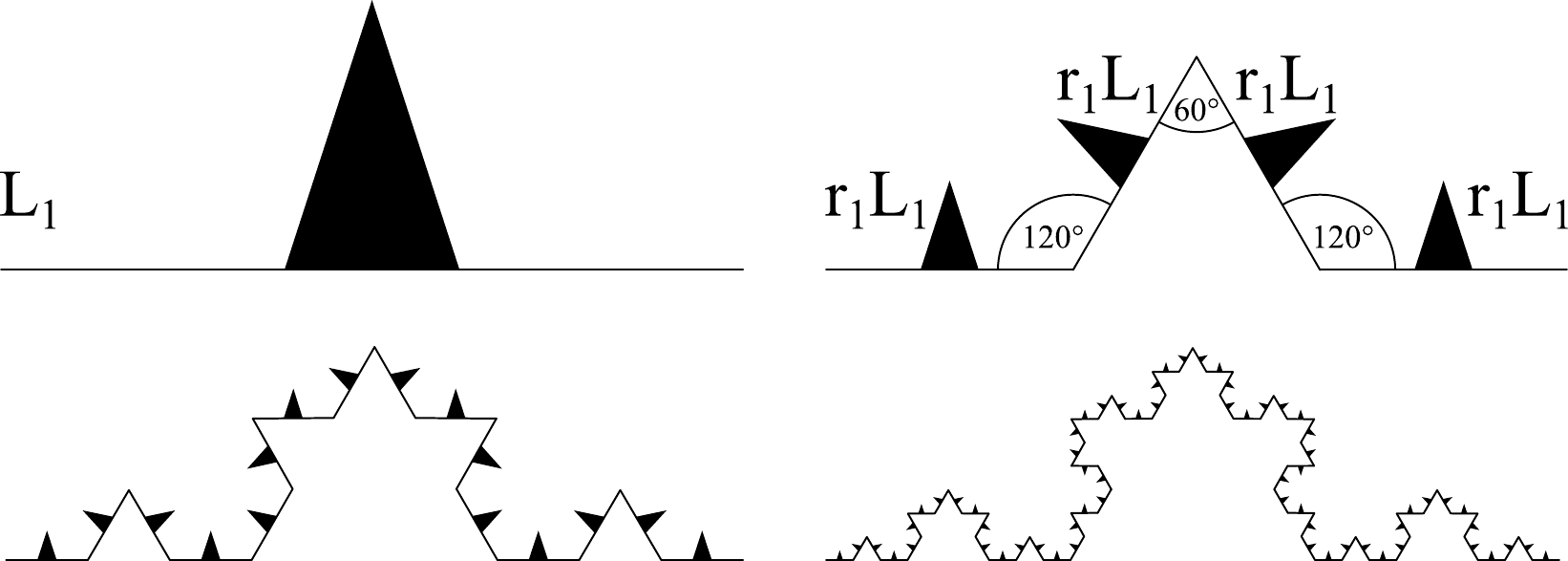}}
\end{center}\caption{\label{Koch}The first three iterations (substitutions) of the Koch
curve as described in \citep{HvK1904,HvK1906,Mandelbrot77a}. The
triangles mark the orientation of the line segments and are not part
of the curve itself.\protect \\
The initiator of the curve is line segment $L_{1}$. The scaling factor
is given by $r_{1}=\nicefrac{1}{3}$. All angles are multiples of
$60\text{\textdegree=}\nicefrac{\pi}{3}rad$. The similarity dimension
is $D=\nicefrac{\log\left(4\right)}{\log\left(3\right)}\simeq1.26186$.}
\end{figure}

Where $n\in\mathbb{N}$ is the number of copies, $0<r<1;\,r\in\mathbb{R}$
is the scaling factor and $D>0;\,D\in\mathbb{R}$ is the similarity
dimension. This equation works for many fractal curves such as the
Koch curve as shown in Figure~\ref{Koch}. Its generator is defined
by a polygonal chain of four line segments of unit length scaled down
by $r_{1}=\nicefrac{1}{3}$.

B. B. Mandelbrot proposed in \citep{Mandelbrot77a} a generalization
for cases where the $n$ line segments of the generator have different
lengths and so different scaling factors $0<r_{i}<1;\,r\in\mathbb{R};\,1\leq i\leq n$.

\begin{equation}
1=\sum_{i=1}^{n}r_{i}^{D}
\end{equation}

\begin{figure}[H]
\begin{center}
\resizebox{0.85\textwidth}{!}{%

\includegraphics{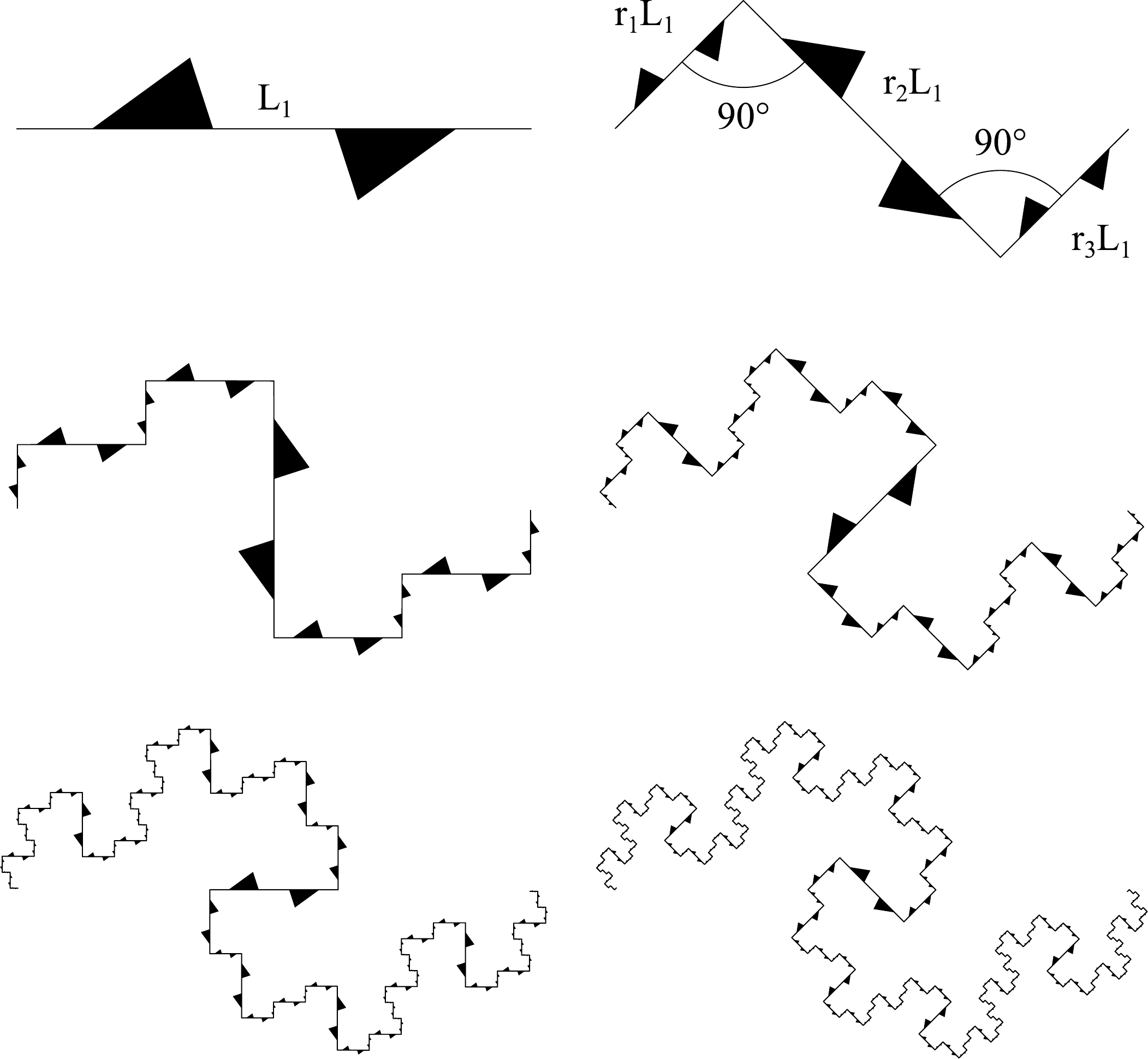}}
\end{center}\caption{\label{Galaxy}The first five iterations (substitutions) of a fractal
curve based on one generator with 3 line segments with two different
lengths. The triangles mark the orientation of the line segments and
are not part of the curve itself.\protect \\
The initiator of the curve is line segment $L_{1}$ as shown on the
upper left.\protect \\
The generator which is equivalent to the first iteration is shown
in the upper right.\protect \\
The scaling factors are given by $2r_{1}=r_{2}=2r_{3}=\nicefrac{1}{\sqrt{2}}$.
All angles are right angles. The similarity dimension is $D\simeq1.52361$.}
\end{figure}

A corresponding example is shown in Figure~\ref{Galaxy}.

This equation can also be written as:

\begin{equation}
1=n\sum_{i=1}^{n}\frac{1}{n}r_{i}^{D}
\end{equation}

Here $n$ is the factor which increases the number of line segments
in each iteration. The line segments with different scaling factors
$r_{i}$ have all the same relative frequency $\nicefrac{1}{n}$.
To describe the case where the line segments have individual relative
frequencies we replace $n$ with $\lambda_{PF}$ and $\nicefrac{1}{n}$
with $f_{i}$ as introduced in Section~\ref{sec:Perron-Frobenius-Theorem}
and derive the proposed definition for the similarity dimension of
fractal curves with multiple generators:

\begin{equation}
1=\lambda_{PF}\sum_{i=1}^{n}f_{i}r_{i}^{D}\label{eq:SimDimForFracCurWithMulGen}
\end{equation}

\section{The Example\label{sec:The-Example}}

\begin{figure}
\begin{center}
\resizebox{0.95\textwidth}{!}{%

\includegraphics{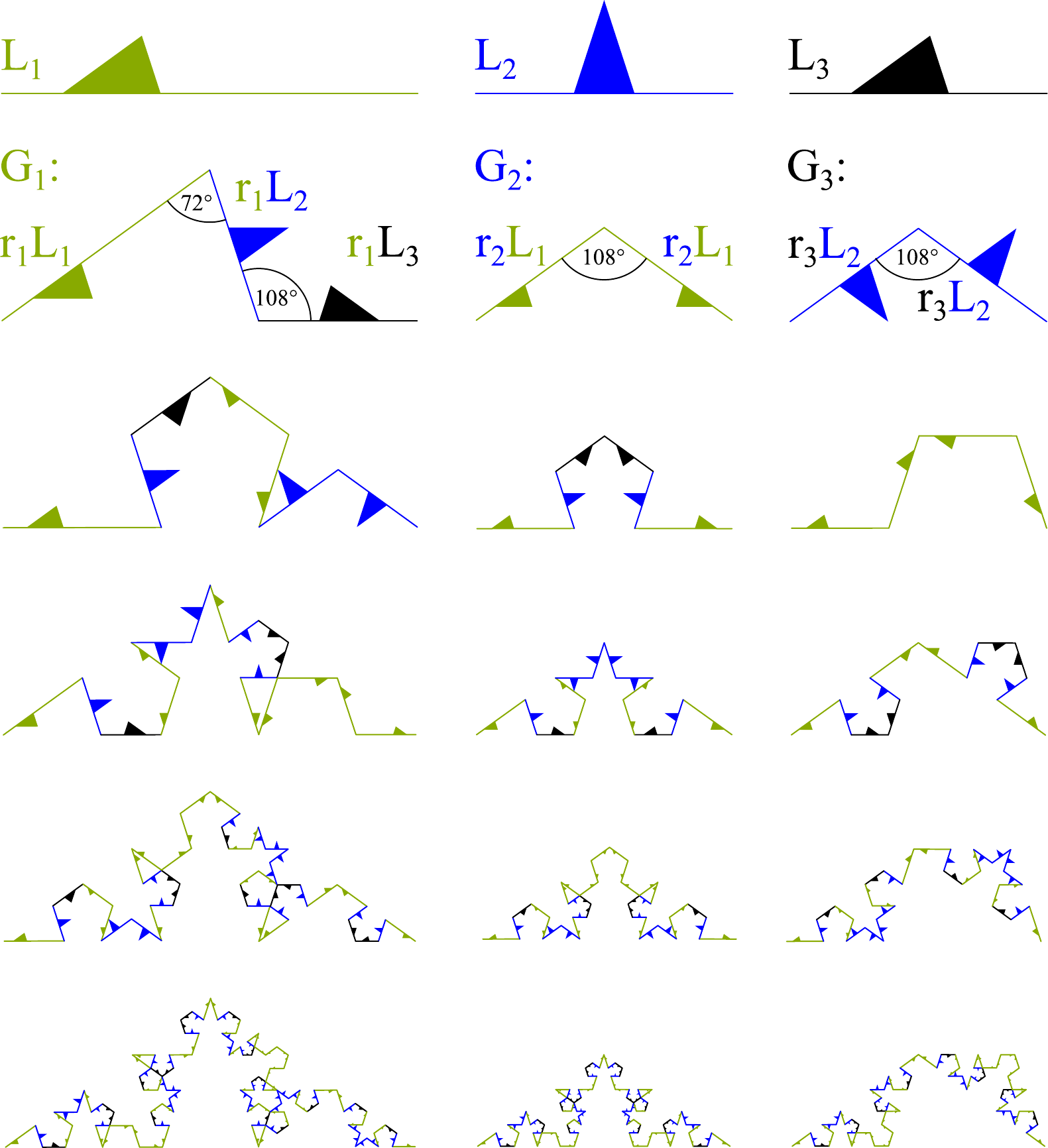}}
\end{center}\caption{\label{generators} The first row shows the three initiators, the
line segments $L_{1}$ (green), $L_{2}$ (blue) and $L_{3}$ (black)
with individual orientations. The triangles mark the orientation of
the line segments and are not part of the curve itself.\protect \\
The second row shows the three corresponding generators $G_{1}$ (green),
$G_{2}$ (blue) and $G_{3}$ (black) which are equivalent to the first
iterations. Each of the generators is defined by a polygonal chain,
made of the line segments $L_{1}$ (green), $L_{2}$ (blue) and $L_{3}$
(black) with orientations, scaled down by $r_{i}$, so that the distance
between start and end point of $G_{i}$ is equal to the length of
the assigned line segment $L_{i}$.\protect \\
The other rows show the result of repeated iterations hence substitutions.
\protect \\
The relative lengths of line segments $L_{1}$, $L_{2}$ and $L_{3}$
are given by $\nicefrac{l_{1}}{l_{2}}=\nicefrac{l_{1}}{l_{3}}=\varphi=\nicefrac{\sin\left(\nicefrac{2\pi}{5}\right)}{\sin\left(\nicefrac{\pi}{5}\right)}\simeq1,618034$.\protect \\
The scaling factors $r_{1}$, $r_{2}$ and $r{}_{3}$ are given by
$r_{1}=r_{3}=\nicefrac{1}{\varphi}=\varphi-1\simeq0,618034$ and $r_{2}=r_{1}^{2}=\nicefrac{1}{\varphi^{2}}=2-\varphi\simeq0,381966$.\protect \\
All angles are multiples of $36\text{\textdegree=}\nicefrac{\pi}{5}\,rad$.\protect \\
The similarity dimension is $D\simeq1.47814$ as derived in Section~\ref{sec:The-Example}.}
\end{figure}

In the following part we discuss an example, in detail the curve shown
in Figure~\ref{curve} and Figure~\ref{generators}. The curve is
defined by a set of three generators or substitution rules $G_{1}$,
$G_{2}$ and $G_{3}$ which are assigned to a set of three lines $L_{1}$,
$L_{2}$ and $L_{3}$ with lengths $l_{1}$, $l_{2}$ and $l_{3}$.
The generators $G_{i}$ describe how a line segment $L_{i}$ can be
replaced by a polygonal chain made of $k_{i,1}$, $k_{i,2}$ and $k_{i,3}$
line segments $L_{1}$, $L_{2}$ and $L_{3}$ scaled down by $r_{i}$.
As a consequence each generator $G_{i}$ and each type of line $L_{i}$
is also assigned to a corresponding inflation multiplier $r_{i}$.

The curve as shown in Figure~\ref{curve} and Figure~\ref{generators}
has the following substitution matrix:

\begin{equation}
M=\left(\begin{array}{ccc}
k_{1,1} & k_{1,2} & k_{1,3}\\
k_{2,1} & k_{2,2} & k_{2,3}\\
k_{3,1} & k_{3,2} & k_{3,3}
\end{array}\right)=\left(\begin{array}{ccc}
1 & 2 & 0\\
1 & 0 & 2\\
1 & 0 & 0
\end{array}\right)
\end{equation}

The columns correspond to the generators $G_{1}$, $G_{2}$ and $G_{3}$
while the rows are assigned to the number of line segments $L_{1}$,
$L_{2}$ and $L_{3}$.\\
Obviously the substitution matrix $M$ is primitive because $M^{3}$
has only positive entries greater zero:

\begin{equation}
M^{3}=\left(\begin{array}{ccc}
9 & 6 & 4\\
5 & 6 & 4\\
3 & 2 & 4
\end{array}\right)>0
\end{equation}

The normalized right Perron-Frobenius eigenvector $x_{PFR}^{*}$ and
the Perron-Frobenius eigenvalue $\lambda_{PF}$ are calculated with
numerical methods:

\begin{equation}
\lambda_{PF}\simeq2.46750
\end{equation}
\begin{equation}
x_{PFR}^{*}=\left(\begin{array}{c}
f_{1}\\
f_{2}\\
f_{3}
\end{array}\right)\simeq\left(\begin{array}{c}
0.46750\\
0.34303\\
0.18946
\end{array}\right)
\end{equation}

Together with Equation~\eqref{eq:SimDimForFracCurWithMulGen} and
the scaling factors $r_{1},r_{2},r_{3}$ in Figure~\ref{generators}
we have:

\begin{equation}
1\simeq2.46750\left(0.46750\left(0,618034\right)^{D}+0.34303\left(0,381966\right)^{D}+0.18946\left(0,618034\right)^{D}\right)
\end{equation}

The similarity dimension $D$ of the curve in Figure~\ref{curve}
and Figure~\ref{generators} is calculated with a numerical method:

\begin{equation}
D\simeq1.47814
\end{equation}

\section{Summary}

Based on B. B. Mandelbrot's equation for the similarity dimension
of fractal curves with one generator and line segments with different
lengths and scaling factors and the tool sets for the analysis of
aperiodic substitution tilings such as the Perron Frobenius theorem,
we derived a proposal for a definition for the similarity dimension
of fractal curves with multiple generators.

\section*{Acknowledgment}

The author would like to thank Dirk Frettlöh (Bielefeld University),
Christian Mayr (Technische Universität Dresden), Klaus-Peter Nischke,
Asta Richter (Technical University of Applied Sciences Wildau), Christoph
Richard (Friedrich-Alexander-Universität Erlangen-Nürnberg), Christian
Richter (Friedrich-Schiller-Universität Jena) and Jeffrey Ventrella
for the helpful discussions and comments.

\bibliographystyle{alpha}
\bibliography{../Paper/girih7-reference___}

\end{document}